\documentclass{amsart}

\def\Q{{\mathbf Q}}
\def\Z{{\mathbf Z}}
\def\C{{\mathbf C}}
\def\r{{\mathfrak G}}
\def\Gal{\mathrm{Gal}}
\def\End{\mathrm{End}}
\def\Aut{\mathrm{Aut}}

\def\fchar{\mathrm{char}}

\def\X{{A}}
\def\Y{{B}}
\def\muA{\mu_\X}

\newtheorem{thm}{Theorem}[section]
\newtheorem{lem}[thm]{Lemma}

\newtheorem{conj}[thm]{Conjecture}

\theoremstyle{definition}
\newtheorem{defn}[thm]{Definition}

\title[Images of $\ell$-adic representations and automorphisms]
{Images of $\ell$-adic representations and automorphisms of
abelian varieties}

\author[A.\ Silverberg]{A.\ Silverberg}
\address{Department of Mathematics, Ohio State University,
231 W.\ 18 Avenue,
Columbus, Ohio 43210--1174, USA}
\email{silver\char`\@math.ohio-state.edu}

\author[Yu. G. Zarhin]{Yu. G. Zarhin}
\address{Department of Mathematics, Pennsylvania State University,
University Park, PA 16802, USA,
\newline
\indent Institute for Mathematical Problems in Biology,
Russian Academy of Sciences, Push\-chino, Moscow Region, 142292, Russia}
\email{zarhin\char`\@math.psu.edu}

\begin{document}

\maketitle

\section{Introduction}

Suppose that $F$ is either a global field or a
finitely generated extension of $\Q$, $\X$ is an abelian variety
over $F$, $\ell$ is a  prime number, and $\ell \neq \fchar(F)$.
Let $\r_\ell(F,\X)$ denote the algebraic envelope of the image
of the absolute Galois group of $F$ under the $\ell$-adic
representation associated to $\X$, and
let $\r_\ell(F,\X)^0$ denote its identity connected component.
In \S\ref{exclasses} we prove that
the intersection of
$\r_\ell(F,\X)^0(\Q_\ell)$ with the torsion subgroup of the center of
$\End(\X) \otimes \Q$ is independent of $\ell$.
In the case where $F$ is a finitely generated extension of $\Q$,
this would follow from the Mumford-Tate Conjecture.
Our results do not assume the Mumford-Tate Conjecture,
and apply even in the positive
characteristic case, where there is no analogue of the
Mumford-Tate Conjecture.
The result in the characteristic zero case can therefore be viewed
as providing evidence in the direction of the Mumford-Tate conjecture.

Let $F_{\Phi,\ell}(\X)$ be the smallest extension $F'$ of $F$ such that
$\r_\ell(F',\X)$ is connected. We call this extension the
$\ell$-connectedness extension, or connectedness extension.

The algebraic group $\r_\ell(F,\X)$
and the field $F_{\Phi,\ell}(\X)$ were
introduced by Serre (\cite{serre},
\cite{resume}, \cite{serremotives}), who
proved that if $F$ is a global field or a finitely
generated extension of $\Q$,
then $F_{\Phi,\ell}(\X)$ is independent of $\ell$
(see \cite{serre},
\cite{resume}, Proposition 1.1 of \cite{LPmathann}, and
Proposition 6.14 of \cite{LPinv}). In such cases,
we will denote the field $F_{\Phi,\ell}(\X)$ by $F_\Phi(\X)$.
Our results rely heavily on the $\ell$-independence results
of Serre and their generalizations \cite{LPinv}.

Let $F(\End(\X))$
denote the smallest extension of $F$ over which all the endomorphisms
of $\X$ are defined. Then
(see Proposition 2.10 of \cite{conn}),
$$F(\End(\X)) \subseteq F_{\Phi,\ell}(\X).$$
By enlarging the ground field, we may assume that
$F = F(\End(\X)) = F_{\Phi,\ell}(\X)$. We then consider twists of
the $\ell$-adic representations associated to $\X$.
The results of this paper follow from the
$\ell$-independence of the connectedness extensions associated
to these twists.

See \cite{connweil} for a study of the connectedness extensions
attached to twists of abelian varieties. See \cite{conn} for
conditions for the connectedness of $\r_\ell(F,\X)$.

\section{Definitions, notation, and lemmas}
\label{notation}

Let $\Z$, $\Q$, and $\C$ denote respectively the integers, rational
numbers, and complex numbers.
If $G$ is an algebraic
group, let $G^0$ denote the identity connected component.
If $F$ is a field, let $F^s$ denote a separable closure of $F$ and let
${\bar F}$ denote an algebraic closure of $F$. If $\X$ is an abelian variety
over a field $F$, write $\End_F(\X)$
for the endomorphisms of $\X$ which are defined over $F$,
write $\Aut_F(\X)$ for the automorphisms of $\X$ defined over $F$, let
$\End(\X) = \End_{F^s}(\X)$, let $\End^0(\X) = \End(\X) \otimes_\Z \Q$,
and let $\End^0_F(\X) = \End_F(\X) \otimes_\Z \Q$.
If $\ell$ is a  prime number and $\ell \neq \fchar(F)$, let
$T_\ell(\X) = {\displaystyle \lim_\leftarrow \X_{\ell^r}}$
(the Tate module), let
$V_\ell(\X) = T_\ell(\X) \otimes_{\Z_\ell}\Q_\ell$, and let
$\rho_{\X,\ell}$
denote the $\ell$-adic representation
$$\rho_{\X,\ell} : \Gal(F^s/F) \to \Aut(T_\ell(\X)) \subseteq
\Aut(V_\ell(\X)).$$
Let $\r_\ell(F,\X)$
denote the algebraic envelope of the image of $\rho_{\X,\ell}$,
i.e., the Zariski closure in $\Aut(V_\ell(\X))$ of the image of
$\Gal(F^s/L)$ under $\rho_{\X,\ell}$.
Let $F_{\Phi,\ell}(\X)$ be the smallest extension $F'$ of $F$
in $F^s$ such that $\r_\ell(F',\X)$ is connected.

\begin{lem}[Lemma 2.7 of \cite{conn}]
\label{conncomp}
If $\X$ is an abelian variety over a field $F$, $L$ is a finite
extension of $F$ in $F^s$, and $\ell$ is a prime number, then
$$\r_\ell(L,\X) \subseteq \r_\ell(F,\X)
\text{ and } \r_\ell(L,\X)^0 = \r_\ell(F,\X)^0.$$
In particular, if $\r_\ell(F,\X)$
is connected, then $\r_\ell(F,\X) = \r_\ell(L,\X)$.
\end{lem}

\begin{lem}
\label{conncomplem}
Suppose $\X$ and $\Y$ are abelian varieties over a field $F$, $L$ is a finite
extension of $F$ in $F^s$, $\ell$ is a prime number,
$\ell \neq \fchar(F)$, $\r_\ell(F,\X)$ is
connected, and $\X$ and $\Y$ are isomorphic over $L$. Then:
\begin{enumerate}
\item[{(i)}] $\r_\ell(F,\Y)^0 = \r_\ell(F,\X)$, and
\item[{(ii)}]$\r_\ell(L,\Y)$ is connected, i.e.,
$F_{\Phi,\ell}(\Y) \subseteq L$.
\end{enumerate}
\end{lem}

\begin{proof}
Since $\X$ and $\Y$ are isomorphic over $L$, and $\r_\ell(F,\X)$ is
connected, we have
$$\r_\ell(L,\Y) = \r_\ell(L,\X) = \r_\ell(F,\X) = \r_\ell(F,\X)^0$$
$$= \r_\ell(L,\X)^0 = \r_\ell(L,\Y)^0 = \r_\ell(F,\Y)^0,$$
using Lemma \ref{conncomp}. The result follows.
\end{proof}

\begin{lem}
\label{strcompat}
Suppose $\X$ is an abelian variety over a global field $F$ and
$$c : \Gal(F^s/F) \to \End^0_F(\X)^\times$$ is a continuous
homomorphism of finite order. For each
prime $\ell \ne \fchar(F)$ let
$$\rho_{\ell,c}: \Gal(F^s/F)  \to \Aut(V_{\ell}(\X)),\quad
 \sigma \mapsto c(\sigma) \rho_{\X,\ell}(\sigma)$$
be the twist of $\rho_{\X,\ell}$. Then $\{\rho_{\ell,c}\}$
constitutes a strictly compatible system of integral
$\ell$-adic representations of $\Gal(F^s/F)$.
More precisely, suppose $M$ is a finite Galois extension of
$F$ such that $c$ factors through $\Gal(M/F)$.
Let $S_\ell$ be the set of finite places $v$ of $F$ such that
either $\X$ has bad reduction at $v$, $v$ is ramified in
$M/F$, or the residue characteristic of $v$ is $\ell$.
Let $v$ be a finite place of $F$,
$w$ a place of $F^s$ lying over $v$, and $\kappa_v$ and
$\kappa_w$ the residue fields at $v$ and $w$, respectively.
Let $\tau \in \Gal(F^s/F)$ be an element that acts as
the Frobenius automorphism of $\kappa_w/\kappa_v$.
Suppose that $v \notin S_\ell$, and let
$$\varphi_w = \rho_{\ell,c}(\tau) \in \Aut(V_{\ell}(\X)).$$
Then $\rho_{\ell,c}$ is unramified at $v$,
the characteristic polynomial
$P_v(t)$ of $\varphi_w$ lies in $\Z[t]$ and does not depend on
the choice of $w$ and $\ell$,
and the roots of $P_v(t)$ all have complex absolute value
$\sqrt{\#\kappa_v}$.
\end{lem}

\begin{proof}
We use that $\{\rho_{\X,\ell}\}$ is a strictly compatible system
of integral $\ell$-adic representations of $\Gal(F^s/F)$
(see \cite{Weil} and I.2 of \cite{abreps}).
Let $\X_v$ be the abelian variety over $\kappa_v$
which is the reduction of $\X$ at $v$.
The choice of $w$ allows us to identify
the Tate modules $V_\ell(\X)$ and $V_\ell(\X_v)$, and this
identification is compatible with the natural embedding
$\End^0(\X) \hookrightarrow \End^0(\X_v)$.
Let
$$Fr_w = \rho_{\X,\ell}(\tau) \in \Aut(V_{\ell}(\X)).$$
Then
$$\varphi_w = c(\tau)Fr_w \in \Aut(V_\ell(\X))
\subseteq \Aut(V_\ell(\X_v)),$$
and the identification of $\Aut(V_\ell(\X))$ with $\Aut(V_\ell(\X_v))$
identifies $Fr_w$ with the Frobenius endomorphism of
$\X_v$ inside $\Aut(V_\ell(\X_v))$.
It follows from Weil's results on endomorphisms of abelian varieties
that $P_v(t)$ has rational coefficients
which do not depend on the choice of $w$ and $\ell$.
If $m = [M:F]$ then
$(c(\tau)Fr_w)^m = (Fr_w)^m \in \End(\X_v)$
and therefore all roots of $P_v(t)$ are algebraic
integers. Therefore, $P_v(t) \in \Z[t]$. Further,
Weil's results imply that the eigenvalues of $\varphi_w$
have absolute value $\sqrt{\#\kappa_v}$.
\end{proof}

\begin{thm}
\label{galois}
If $F$ is either a finitely generated extension of $\Q$ or
a function field in one variable over a finite field,
then every finite abelian group occurs as
a Galois group over $F$.
\end{thm}

\begin{proof}
See Theorem 3.12c of \cite{Saltman}, and
IV.2.1 and IV.1.2 of \cite{Matzat}.
\end{proof}

Next we define the Mumford-Tate group of a complex abelian variety $\X$
(see \S2 of \cite{Ribet} or \S6 of \cite{Izv}). If $\X$ is a complex abelian
variety, let $V = H_1(\X(\C),\Q)$ and consider the Hodge decomposition
$V \otimes \C = H_1(\X(\C),\C) = H^{-1,0} \oplus H^{0,-1}$.
Define a homomorphism $\mu : {\mathbf G}_m \to GL(V)$ as follows. For
$z \in \C$, let $\mu(z)$ be the automorphism of $V \otimes \C$ which is
multiplication by $z$ on $H^{-1,0}$ and is the identity on $H^{0,-1}$.

\begin{defn}
The {\em Mumford-Tate group} $MT_\X$ of $\X$ is the smallest
algebraic subgroup of $GL(V)$, defined over $\Q$, which after extension of
scalars to $\C$ contains the image of $\mu$.
\end{defn}

If $\X$ is an abelian variety over a subfield $F$ of $\C$, we fix an
embedding of ${\bar F}$ in $\C$. This gives an identification of
$V_\ell(\X)$ with $H_1(\X,\Q)\otimes\Q_\ell$, and allows us to view
$MT_\X \times \Q_\ell$ as a linear $\Q_\ell$-algebraic subgroup of
$GL(V_\ell(\X))$. Let
$$MT_{\X,\ell} = MT_\X \times_\Q \Q_\ell.$$
Then $MT_{\X}(\Q_\ell) = MT_{\X,\ell}(\Q_\ell)$.
The Mumford-Tate conjecture for abelian varieties (see \cite{serrereps})
may be reformulated as the equality of $\Q_\ell$-algebraic groups,
$\r_\ell(F,\X)^0 = MT_{\X,\ell}$.

\begin{conj}[Mumford-Tate Conjecture]
\label{mtconj}
If $\X$ is an abelian variety over a finitely generated extension $F$ of
$\Q$, then $\r_\ell(F,\X)^0 = MT_{\X,\ell}$.
\end{conj}

The inclusion $\r_\ell(F,\X)^0 \subseteq MT_{\X,\ell}$ was
proved by Piatetski-Shapiro \cite{ps}, Deligne \cite{sln900},
and Borovoi \cite{Borovoisb}.

It is well-known that $MT_A$ contains the homotheties
${\mathbf G}_m$ and that the centralizer of
$MT_A$ in $\End(V)$ is $\End^0(A)$. Therefore, the center of
$MT_A(\Q)$ contains $-1$ and is contained in the center of $\End^0(A)$.

\section{$\ell$-independence}
\label{exclasses}

Suppose that $F$ is either a
finitely generated extension of $\Q$ or a global field.
Suppose $F = F_\Phi(\X)$, so that
$\r_\ell(F,\X) = \r_\ell(F,\X)^0 = \r_\ell(L,\X)$
for all finite extensions $L$ of $F$.
It follows from \cite{faltings}, \cite{fw}, \cite{z0}, \cite{z1},
and VI.5 and XII.2 of \cite{Mori2}
that $\r_\ell(F,\X)$ is a reductive
$\Q_\ell$-algebraic group, whose centralizer in $\End(V_\ell(\X))$
is $\End(\X) \otimes \Q_\ell$. This implies that the center
of $\r_\ell(F,\X)(\Q_\ell)$ is contained in
$(Z \otimes \Q_\ell)^\times$, where $Z$ is the center
of $\End^0(\X)$.
Let $\muA$ denote the group of elements of finite order in
the center of $\End^0(\X)$.
In the case where $F \subset \C$,
the Mumford-Tate Conjecture (Conjecture \ref{mtconj})
would imply that
$\muA \cap \r_\ell(F,\X)(\Q_\ell)$
is the torsion subgroup of the center of $MT_\X(\Q)$,
and therefore is independent of $\ell$.
In the following two results we prove that
$\muA \cap \r_\ell(F,\X)(\Q_\ell)$
is independent of $\ell$ (without assuming the
Mumford-Tate Conjecture).

It follows from Weil's results on abelian varieties \cite{Weil}
(as was pointed out by Deligne; see 2.3 of \cite{serrereps}) that
$\r_\ell(F,\X)$ contains the homotheties ${\mathbf G}_m$.
In particular,
$$-1 \in \r_\ell(F,\X)(\Q_{\ell}).$$

\begin{thm}
\label{discondcor}
Suppose $\X$ is an abelian variety over a finitely generated
extension $F$ of $\Q$,
and $F = F_{\Phi}(\X)$. Let $\muA$ denote the
group of elements of finite order in
the center of $\End^0(\X)$.
Then $\muA \cap \r_\ell(F,\X)(\Q_\ell)$ is independent of
the prime $\ell$.
\end{thm}

\begin{proof}
Over $\C$, we can view $A$ as $\C^d/L$ with $L$ a lattice in
$\C^d$. Then $L'=\sum_{\gamma \in \muA}\gamma(L)$ is a
$\muA$-invariant lattice in $\C^d$ that contains $L$ as a
subgroup of finite index.
The complex abelian variety $\C^d/L'$
has a model $A'$ defined over a finite extension
$F'$ of $F$ such that $\X$ and $A'$ are $F'$-isogenous and
$\muA$ coincides with the set of elements of finite order in
the center of $\End(A')$. Since
$\r_\ell(F',A') = \r_\ell(F',\X) = \r_\ell(F,\X)$,
we may assume without loss of generality that $\muA$ coincides
with the set of elements of finite order in
the center of $\End(\X)$.
By Theorem \ref{galois},
we can choose an abelian extension $M$ of $F$ such that
$\Gal(M/F)$ is isomorphic to $\muA$. Let
$$\chi : \Gal(M/F) \to \muA$$
be an isomorphism, let $c : \Gal(F^s/F) \to \muA$ be the
composition of $\chi$ with the projection
$\Gal(F^s/F) \to \Gal(M/F)$, and
let $\Y$ denote the twist of $\X$ by the cocycle
induced by $c$. By Lemma \ref{conncomplem}i,
$\r_\ell(F,\Y)^0 = \r_\ell(F,\X)$.
The character $c$ induces an isomorphism
$$\Gal(M/F_{\Phi}(\Y)) \cong \muA \cap \r_\ell(F,\Y)^0(\Q_\ell) =
\muA \cap \r_\ell(F,\X)(\Q_\ell).$$
Since $\Gal(M/F_{\Phi}(\Y))$ is independent of $\ell$, we are done.
\end{proof}

\begin{thm}
\label{discondcor2}
Suppose $F$ is a function field in one variable over
a finite field, $\X$ is an abelian variety over $F$, and
$\ell$ is a prime number not equal to $\fchar(F)$.
Suppose $F = F_{\Phi}(\X)$, and let $\muA$ denote the
group of elements of finite order in the center of $\End^0(\X)$.
Then $\muA \cap \r_\ell(F,\X)(\Q_\ell)$ is independent of $\ell$.
\end{thm}

\begin{proof}
By Theorem \ref{galois} we
can choose an abelian extension $M$ of $F$ such that
$\Gal(M/F)$ is isomorphic to $\muA$. Let
$$\chi : \Gal(M/F) \to \muA$$
be an isomorphism, let $c : \Gal(F^s/F) \to \muA$ be the
composition of $\chi$ with the projection
$\Gal(F^s/F) \to \Gal(M/F)$, and
define $\rho_{\ell,c} : \Gal(F^s/F) \to \Aut(V_\ell(A))$
by $\rho_{\ell,c}(\sigma) = c(\sigma)\rho_{A,\ell}(\sigma)$.
For $F \subseteq F' \subseteq F^s$, let $\r_{\ell,c}(F')$ denote
the Zariski closure of
$\rho_{\ell,c}(\Gal(F^s/F'))$. Let $F_{\Phi,c}$
denote the smallest extension $F'$ of $F$ in $F^s$ such that
$\r_{\ell,c}(F')$ is connected. By Lemma \ref{strcompat},
$\{\rho_{\ell,c}\}$ is a strictly compatible system of
integral $\ell$-adic representations.
Therefore by Proposition 6.14 of \cite{LPinv},
$F_{\Phi,c}$ is independent of $\ell$.
By definition,
$$\Gal(M/F_{\Phi,c}) \cong
\muA \cap \r_{\ell,c}(F)^0(\Q_\ell) \quad {\text{ and }}
\quad \r_{\ell,c}(M) = \r_\ell(M,A).$$
Lemma \ref{conncomp} is valid with $\r_{\ell,c}(F')$ in place of
$\r_\ell(F',A)$; the proof remains unchanged. Therefore,
$$\r_{\ell,c}(F)^0 = \r_{\ell,c}(M)^0 = \r_\ell(M,A)^0 = \r_\ell(F,A).$$
Since $\Gal(M/F_{\Phi,c})$ is independent of $\ell$, we are done.
\end{proof}

\end{document}